\documentclass[11pt]{article}
\usepackage{amsfonts}
\usepackage{amsfonts} 
\usepackage{pdfsync} 

\usepackage{amsthm}
\usepackage{amsmath}
\usepackage{color}

\title{Mixed boundary value problems for fully nonlinear degenerate or singular equations} 

\author{Isabeau Birindelli \\
Dipartimento di Matematica, Sapienza Universit\`a\  di Roma
\and
  Fran\c{c}oise Demengel\\
  D\'epartement de Math\'ematiques,
Universit\'e\  de Cergy-Pontoise
  \and
   Fabiana  Leoni\\ 
   Dipartimento di Matematica, Sapienza Universit\`a\  di Roma}

\date{}

\catcode`@=11 \@addtoreset{equation}{section} \catcode`@=12

\newtheorem{theo}{Theorem}[section]
\newtheorem{prop}[theo]{Proposition}
\newtheorem{rema}[theo]{Remark}
\newtheorem{defi}[theo]{Definition}
\newtheorem{cor}[theo]{Corollary}
\newtheorem{lemme}[theo]{Lemma}

\def\R{\mathbb  R}

\begin{document}
\maketitle
\begin{abstract}
We prove existence, uniqueness and regularity results for mixed boundary value problems associated with fully nonlinear, possibly singular or degenerate elliptic equations. Our main result is a global H\"older estimate for solutions, obtained by means of the comparison principle and the construction of ad hoc barriers. The global H\"older estimate immediately yields a compactness result in the space of solutions, which could be applied in  the study of principal eigenvalues and principal eigenfunctions of mixed boundary value problems.
\smallskip

\end{abstract}

\section{Introduction}
In this article we prove existence, uniqueness  and regularity of solutions to a mixed boundary value problem for singular or degenerate  fully nonlinear equations. 
The model problem we have in mind is the following: let   $\Omega$ denote a  bounded domain in $\R^N$ and let ${\cal N} $ be some  relatively open subset, of class $C^2$, of $\partial \Omega$. We  set $\mathcal{D}= \partial \Omega \setminus  {\mathcal N}$ and consider 
  \begin{equation}\label{exmix} 
 \left\{ \begin{array}{c} 
-| \nabla u|^\alpha F( x, \nabla u, D^2 u) +\beta(u)=  f \qquad  {\rm in} \ \Omega,  \\
 u = 0  \ {\rm on }\  \mathcal{D},\quad  \displaystyle \frac{\partial u}{\partial \vec n} = 0 \ { \rm on} \ {\mathcal N}
 \end{array}\right.
\end{equation} 
with $\alpha>-1$, $F$ uniformly elliptic and $\beta:\R\to \R$ continuous and increasing. The typical example of operator $F$ we have in mind is one of  Pucci's extremal operators, defined as
 $$
{\mathcal M}^+_{a,A} (M)= A\, \sum_{i=1}^N \lambda_i^+ - a\, \sum_{i=1}^N\lambda_i^-
$$
and
$$
{\mathcal M}^-_{a,A} (M)= a\, \sum_{i=1}^N \lambda_i^+ - A\, \sum_{i=1}^N\lambda_i^-\, ,
$$
where $A\geq a>0$, $\lambda_1, \ldots ,\lambda_N$ are the eigenvalues of the matrix $M$, and $\lambda_i^+=\max\{\lambda_i,0\}$, $\lambda_i^-=\max\{ -\lambda_i,0\}$. 
Operators ${\mathcal M}_{a,A}^{\pm}$ act both as prototypes and barriers for the class of uniformly elliptic operators $F$ having ellipticity constants $A$ and $a$. 
We will further suppose that $f$ is bounded and continuous.

Mixed boundary conditions come up naturally. A very basic but enlightening  example is that of a drum whose batter head may be attached to the hoop through "nails": the Dirichlet part corresponds to the nails and the Neumann part to the rest of the boundary. Another example comes from the necessity to approximate a Neumann boundary problem in unbounded domains $\mathcal O$ with problems posed in bounded domains, that are e.g. intersection with  balls of radius $R$ i.e. $\mathcal O\cap B_R$. Then, it is natural to impose Dirichlet on the boundary of the ball and Neumann on $\partial\mathcal O$, see e.g. \cite{Ro}.

The aim is to prove that there exists a unique solution of \eqref{exmix}  and that the solution satisfies a global H\"older estimate. 
Concerning the global H\"older estimate, which is done in Section \ref{section4}, it is important to notice that the difficulty and novelty lies mainly near the points that are 
at the intersection of the Dirichlet part of the boundary with the closure of the Neumann part. In fact, for some of the results we will suppose,  in addition, that 
$\mathcal{D} = \overline{ {\rm int} \mathcal{D}}$, and also that $\mathcal{D}\cap \overline{\cal N}$ is the union of  $(N-2)$-dimensional manifolds of class $C^2$. 
The proof of the H\"older regularity, which is optimal in that case, requires the construction of an ad hoc barrier that is not at all standard, see Lemma \ref{super2}.
Elsewhere the proof is technical, but it follows the lines of other similar results obtained in the Dirichlet or in the Neumann case.

In order to prove existence of solutions via Perron's method \lq\lq \`a la Ishii'', we need to prove a comparison principle and construct sub and super solutions. The comparison principle itself 
requires interior H\"older estimates. These estimates are stated in Subsection \ref{prelreg} but proved in the appendix. We wish to point out   that   Theorem \ref{compprinc}  could be proved  directly, by using the results of \cite{IL}, however  Theorem \ref{lipN},  which provides an estimate on the difference between sub an super solutions,   is interesting in itself and it does not 
assume any sign condition on the supremum of that difference. Hence  it applies to solutions, thus yielding   locally Lipschitz regularity on $\Omega \cup { \cal N}$ for  
\emph{solutions}.  The construction of the sub and super solution is not standard, due to the non linearity and the singular or degenerate character of the operator, jointly with the mixed boundary condition.

 The  existence of solution for analogous problems equipped with Dirichlet  boundary conditions  has been studied in \cite{BD1},  while the Neumann case is treated  when $\alpha = 0$ in \cite{CIL}, and for any $\alpha >-1$ in \cite{Pa}. On  regularity results up to the boundary, which extend the regularity results of   \cite{IS} , let us recall the results of  \cite{BDCocv} for the Dirichlet  case, and those of  \cite{Ri},  \cite{BV} for the Neumann problem.

  Concerning mixed boundary conditions for   fully nonlinear equations, almost nothing is  known. Even in the variational case and even for the Laplacian case,  one can remark that the Neumann condition is not included in the constraint : it is concealed and obtained  a posteriori.  Note also that,  of course, at the intersection of the Neumann and the Dirichlet part, the ${ \cal C}^1$ regularity of solutions cannot hold, since Hopf Principle contradicts this regularity in most of the cases. 
  
In the linear variational setting, we should mention the work of Lieberman \cite{L1,L2,L3} and  the results of Miranda \cite{Mir} and Azzam and Kreyszig \cite{AK}. 
In particular,  in his famous book, Lieberman \cite{L3}  treats the  Laplacian case by  using merely comparison arguments. For the strong maximum principle, we recall the work of \cite{Dav}.  Finally, the case of the $p-$Laplacian, which share some similarities with our case, has been treated in \cite{LLC}, but their procedure is completely variational.

Our interest for the mixed boundary problems came up through a talk of Ireneo Peral, who had considered in \cite{Peral} an extension to fractional Laplacian of the intriguing result of 
Denzler \cite{Denz} on the minimization of the principal eigenvalue with respect to the configuration of the Neumann part and the Dirichlet part of the boundary. 
 At that point, realizing how little was known in a nonlinear non variational setting, it was completely natural to wish to fill the gap. 
The compactness result  we prove at the end of the present paper can be   used to prove the existence of a principal eigenfunction, once we will have defined the principal eigenvalue, on the model of \cite{BNV}. The study of existence and  further properties of  principal eigenfunctions for the mixed problem  will be the subject of a forthcoming paper.

The paper is organized as follows. 
    
     In Section 2  we   precise the notation and the assumptions, and we recall the definition of viscosity solutions adapted to the present context, taking into account the  boundary condition. We will consider viscosity solutions satisfying the Neumann condition in the so called \lq\lq strong sense", and we will show this is equivalent to consider solutions satisfying the Neumann  condition in the generalized viscosity sense. 
          
     In Section 3 we first state  local Lipschitz type estimates between sub and super-solution around points belonging to the Neumann part of the boundary. This will lead to establish  the comparison principle.  We next prove an existence result for the mixed problem, by giving the explicit  construction of sub- and super-solutions satisfying the boundary conditions. 
     
     In Section 4 we provide a global H\"older result,  under an additional regularity assumption on the set given by the intersection between the Dirichlet  and the Neumann part. This assumption  is similar to the so called \lq\lq $\Sigma$ wedge  condition" of Lieberman, \cite{L1}, \cite{L2}. From this we derive  a global compactness result,  as well as an existence result for equations not having zero order terms.   
       
     Finally in Section 5, the Appendix, we give the proof of the local Lipschitz regularity stated in Section 3.

     \medskip
Let us finally  emphasize  the fact that most of the results here enclosed can  be extended to the case of mixed Dirichlet/oblique derivative boundary problems. 

  \section{Notation and  assumptions}

 In all the paper $\alpha$ denotes some real number $>-1$,  $\mathcal{S}$ is the space of symmetric $N\times N$-matrices equipped with the usual ordering   and 
  $\Omega$ denotes an open bounded domain in $\R^N$.  ${\cal N} $ is some  relatively open subset, of class $C^2$, of $\partial \Omega$, and $\mathcal{D}= \partial \Omega \setminus  {\mathcal N}$.  We will suppose   for some of the results in addition that 
  $\mathcal{D} = \overline{ {\rm int} \mathcal{D}}$, and also that $\mathcal{D}\cap \overline{\cal N}$ is the union of  $(N-2)$-dimensional manifolds of class $C^2$. We will denote by $\vec n(x)$ the {\color{red} unit outer}  normal  at any point $x\in \mathcal{N}$, and by $d(x)$ the distance function from $\partial \Omega$. We observe that $-\nabla d(x)=\vec n(x)$ for all $x\in \mathcal{N}$.
\smallskip
 
On the operator
    $F:\overline{\Omega}\times \R^N\times \mathcal{S}\to \R$ we assume that
\begin{itemize}        
 \item[(H1)] 
            $ F( x, p, O)\equiv 0\, ;$
\item[(H2)]   there exist $0< a\leq A$ such  that, for any $x\in \overline{\Omega}$, $p\in \R^N$ and $M, P\in \mathcal{S}$ with $P \geq 0$
            $$a\, {\rm tr}(P) \leq F(x, p, M+ P)-F(x, p, M) \leq A\, {\rm tr}(P)\, ;$$
\item[(H3)]  $F$ is uniformly Lipschitz  continuous in $p$, that is  there exists a constant ${\rm lip}_p F$ such that for any $p, q\in \R^N$ and for any $x\in  \overline{\Omega}$ and $M\in \mathcal{S}$,
             $$|F( x, p, M)-F( x, q, M) | \leq  {\rm lip}_p F |p-q|\, ;$$
 \item[(H4)]  there exist $\theta_F \in ]\frac{1}{2} ,1[$  and a constant $c_F>0$ such that for any $(x,y)\in \overline{ \Omega}^2$ and for all $p\in \R^N$, $M\in \mathcal{S}$
               $$ |F( x, p, M)-F( y, p, M) | \leq c_F |M|\, |x-y|^{\theta_F}\, ; $$
               (Here $|M|$ denotes the  spectral radius of $M$, for example).  
 \end{itemize}
It is easy to see that,   by assumptions (H1) and (H2),  for all $x\in \overline{\Omega},\ p\in \R^N$ and $M\in {\mathcal S}$ one has
\begin{equation}\label{estiPucci}
{\mathcal M}_{a,A}^-(M)\leq F(x,p,M) \leq {\mathcal M}_{a,A}^+(M)\, .
\end{equation}

However, in the results given in the next sections the dependence of operator $F$ on the space variable $x$ and the vector $p$ is allowed, under the assumptions (H3) and (H4) respectively. In this regard, let us observe that assumption (H4) could be replaced by the  hypothesis
\begin{itemize}
\item[(H4)'] 
    for $q>\sup \{ 2, \frac{\alpha + 2}{\alpha+1}\}$, 
               there exists  $\omega:[0,+\infty)\to [0,+\infty)$, satisfying $\omega(0)=0$ and being continuous at $0$, such that for any $(x_j, y_j) \in (\overline{ \Omega})^2$ and $X_j, Y_j$ satisfying 
              $$\left( \begin{array}{cc} 
              X_j& 0\\
              0& Y_j  \end{array}\right) \leq 3j |x_j-y_j|^{ q-2} \left(\begin{array}{cc}  I& -I\\
              -I& I\end{array} \right)\, ,$$
one has
     $$F (x_j,  j |x_j-y_j|^{ q-2} (x_j-y_j), X_j) -F(y_j, j |x_j-y_j|^{ q-2} (x_j-y_j), -Y_j) \leq \omega ( j |x_j-y_j|^q).$$
\end{itemize}
Assumption (H4)' was introduced in \cite{CIL}, see also \cite{BD1}, and it is the natural assumption which makes  the standard proof of the comparison   principle for viscosity solutions works. 
\medskip
                                                                                                                                             
  Next, we precise what we mean by a sub- or a super-solution of the mixed boundary value problem 
  \begin{equation}\label{defmix} 
 \left\{ \begin{array}{c} 
-| \nabla u|^\alpha F( x, \nabla u, D^2 u) + \beta(u)=  f \qquad  {\rm in} \ \Omega,  \\[1ex]
 u = 0  \ {\rm on }\  \mathcal{D},\quad  \displaystyle \frac{\partial u}{\partial \vec n} = 0 \ { \rm on} \ {\mathcal N}
 \end{array}\right.
\end{equation} 
where $\beta:\R \to \R$ is a nondecreasing continuous function, which will be always assumed to satisfy $\beta(0)=0$ without loss of generality, and the datum $f:\overline{\Omega}\to \R$ is a bounded function.

We denote by $USC (\overline \Omega)$ and $LSC(\overline \Omega)$ the set of respectively  upper and lower semicontinuous functions on the set $\overline \Omega$.
 \begin{defi}\label{defvis}
 A function $u\in USC(\overline \Omega)$   is a sub-solution of  \eqref{defmix} if 
 \begin{itemize}
 \item[{\rm (i)}] $u(x)\leq 0$ for all $x\in \mathcal{D}$;
\item[{\rm (ii)}] for any  $\bar x \in \mathcal{N}$ and  for any function $\varphi$ of class $C^1$ in a neighborhood of $\bar x$ such that $u-\varphi$ has a local maximum point at $\bar x$ one has
$$
\frac{\partial \varphi}{\partial \vec n} (\bar x)\leq 0;
$$
\item[{\rm (iii)}] for any $\bar x \in \Omega$ and for any function $\varphi$ of class $C^2$ in a neighborhood of $\bar x$ such that $u-\varphi$ has a local maximum point at $\bar x$, and such that $\nabla \varphi (\bar x)\neq 0$, one has
 $$-| \nabla \varphi(\bar x) |^\alpha F( \bar x, \nabla \varphi (\bar x), D^2 \varphi( \bar x)) +\beta (u(\bar x)) \leq f (\bar x);$$
 \item[{\rm (iv)}] for any $\bar x \in \Omega$ such that $u$ is locally constant in a neighborhood of $\bar x$ one has
$$ \beta (u(\bar x)) \leq f( \bar x). $$
\end{itemize}                        
   A function $v\in LSC(\overline \Omega)$ is a super-solution if it satisfies the   obvious symmetric  conditions. 
    \end{defi}                                         
 \begin{rema}
 {\rm This \lq\lq weird"  definition  is necessary only for the case  $\alpha <0$, due to the fact that the operator appearing in  the equation is not defined when the gradient is zero. In the case $\alpha \geq 0$, the  classical definition of viscosity solution   is equivalent to  Definition \ref{defvis}, see \cite{AR}.}
 \end{rema}
 
 \begin{rema}\label{weakstrong}
{\rm  Let us emphasize that in Definition \ref{defvis} the Neumann boundary condition on the portion $\mathcal{N}$ is imposed in the the so called  "strong viscosity sense".  However, in the present framework, it is equivalent to the "weak viscosity sense", that is: for any $\bar x\in \mathcal{N}$  and for any   function $\varphi$ of class $C^2$ in a neighborhood of $\bar x$ such that $u-\varphi$ has a local maximum point at $\bar x$, and such that $\nabla \varphi (\bar x)\neq 0$, one has
$$
\min \left\{ \frac{\partial \varphi}{\partial \vec n}(\bar x), -| \nabla \varphi(\bar x) |^\alpha F( \bar x, \nabla \varphi (\bar x), D^2 \varphi( \bar x)) + \beta (u(\bar x))- f (\bar x)\right\} \leq 0\, .
$$
The equivalence of the two conditions can be proved by arguing as  in \cite{Ba1} for the case $\alpha = 0$.  We reproduce    the proof here for the reader's convenience : let us assume that $u$ is a sub-solution in the weak viscosity sense, and let $\bar x\in \mathcal{N}$ and $\varphi (x)$ be such that $u-\varphi$ has a local maximum point at $\bar x$. Assume first that $\varphi$ is of class $C^2$ in a neighborhood of $\bar x$ and suppose, by contradiction, that
$$
\frac{\partial \varphi}{\partial \vec n}(\bar x)>0\, .
$$
Note that this implies in particular that $\nabla \varphi (\bar x)\neq 0$. For $\lambda , \ \mu>0$ to be suitably chosen, let us consider the perturbed test function
$$
\varphi_{\lambda, \mu}(x)= \varphi (x)+\lambda d(x)-\mu\, d^2(x)\, ,
$$
where $d(x)$ denotes the distance function from $\partial \Omega$. Then, $\varphi_{\lambda, \mu}$  is of class $C^2$ in a neighborhood of $\bar x$ and $u- \varphi_{\lambda, \mu}$ still has a maximum point at $\bar x$. Moreover, for $0<\lambda <\frac{\partial \varphi}{\partial \vec n}(\bar x)$, one has
$$
\frac{\partial \varphi_{\lambda, \mu}}{\partial \vec n}(\bar x)= \frac{\partial \varphi}{\partial \vec n}(\bar x)-\lambda >0
$$
and, therefore, by the weak viscosity condition,
$$
\begin{array}{c}
 \displaystyle -|\nabla \varphi(\bar x)-\lambda n(\bar x)|^\alpha F( \bar x, \nabla \varphi (\bar x)-\lambda n(\bar x), D^2 \varphi( \bar x) +\lambda D^2 d(\bar x)-2 \mu n(\bar x)\otimes n(\bar x)) \qquad \qquad   \\[2ex]
\qquad \qquad \qquad \displaystyle +\beta(u(\bar x)) \leq  f (\bar x)\, .
\end{array}
$$
On the other hand, by the uniform ellipticity condition (H2), the previous inequality cannot hold for $\mu$ large enough, and we reach a contradiction showing that $\varphi$ has to satisfy $\frac{\partial \varphi}{\partial \vec n}(\bar x)\leq 0$ (see also \cite{Ba1, Ba2, CIL}). The same conclusion can be established, by density, for $\varphi$ of class $C^1$.
}\end{rema}

  \section{Existence and uniqueness results}
  \subsection{Preliminary regularity results}\label{prelreg}
In this section we state the local regularity result that is needed in order to prove   the comparison result, Theorem \ref{compprinc}.
Let us   point out   that   Theorem \ref{compprinc}  could be proved directly, by using the result of \cite{IL} stating that, 
if $u$ and $v$ are respectively a sub- and a super-solution, and  $\sup_\Omega ( u-v) >0$, then  there exists $C>0$ such that,
 for all $(x,y) \in \overline{ \Omega}^2$,
    $u(x)-v(y) \leq \sup ( u-v) + C\,  |x-y|$, this inequality being a crucial argument in the proof of the comparison principle. 
   However  the theorem below does not assume any sign condition on $\sup (u-v)$ and it is interesting in itself,  since  it allows the case $u = v$, thus yielding   locally Lipschitz regularity on $\Omega \cup { \cal N}$ for  \emph{solutions}.  

The results of the present section are stated for equations not containing zero order terms, since, if present,  they can be moved to the right hand side and treated as part of the forcing terms.

 Let $\bar x\in \mathcal N$ and assume that $B(\bar x,\rho) \cap \partial \Omega \subset {\cal N}$.  
 Suppose that $u$ is  $USC (\overline{B(\bar x,\rho) \cap  \Omega})$ and satisfies 
     $$\left\{ \begin{array}{c}
      -| \nabla u |^\alpha F(x, \nabla u,  D^2 u)  \leq f \  \hbox{in} \  B(\bar x,\rho) \cap  \Omega \\[1ex]
    \displaystyle    \frac{\partial u}{\partial \vec n} \leq 0 \  \hbox{on} \  B(\bar x,\rho) \cap \partial \Omega 
             \end{array}\right.$$
as well as  $v$ is $LSC(\overline{B(\bar x,\rho) \cap  \Omega})$ and satisfies 
       $$\left\{ \begin{array}{c}
      -| \nabla v |^\alpha F(x, \nabla v,  D^2 v)  \geq g \ \hbox{in} \ B(\bar x,\rho) \cap  \Omega\\[2ex]
 \displaystyle \frac{\partial v}{\partial \vec n} \geq 0 \  \hbox{on} \  B(\bar x,\rho) \cap \partial \Omega 
      \end{array}\right.$$
with $f$ and $g$ bounded.
 \begin{theo}\label{lipN}  Under the above conditions for any   $r<\rho$, there exists $c_{r}>0$, depending on $\Omega$, $N$, the ellipticity constants of $F$ and on $\|f\|_\infty,\ \|g\|_\infty$, $\sup u$ and $\inf v$, such that, 
       for all $(x, y)\in (B(\bar x,r) \cap \overline{ \Omega} )^2$, one has
       $$ u(x)-v( y) \leq \sup_{ \overline{ \Omega\cap B(\bar x,\rho)}}  ( u-v) + c_{r} | x-y| \,.$$
\end{theo}       
As it is often the case, the proof of  this Lipschitz regularity will be done once we obtain the analogous H\"older result. This is done in section \ref{appendix}  where we state and prove Proposition \ref{holderN} and then we give the proof of Theorem \ref{lipN}. 

\begin{rema}\label{analogies}
{\rm The estimate provided by Theorem \ref{lipN} is the localized version of the analogous estimate obtained in \cite{Pa} in the case of global Neumann boundary condition. An analogous local estimate can be proved around points belonging to either $\Omega$ or the interior of ${\mathcal D}$. In these latter cases the same proof given for Theorem \ref{lipN} can be applied, with simplified arguments, since the auxiliary function $\psi$ can be chosen of the form
$$
\psi(x,y) = u(x)-v(y)-\sup_{\overline{ \Omega\cap B(\bar x,\rho)}} (u-v) - (M \omega ( | x-y|) + L_1 |x-x_0|^2+ L_1 |y-x_0|^2)
\, ,$$
see \cite{BDCocv} . Let us emphasize that, however,  at this level global H\"older or Lipschitz estimates are still missing, because no result is given around points of ${\mathcal D}\cap \overline{\mathcal N}$. This will be the contribution of Section 4.
}
\end{rema}
\subsection{Comparison Theorem}
We are now in a position to prove the following comparison theorem.
\begin{theo}\label{compprinc}  Let $\Omega\subset \R^N$ be a bounded domain satisfying  the uniform exterior sphere condition and let  ${\cal N}\subset \partial \Omega$ be a 
relatively open part  of class ${ \cal C}^2$. Let further    $\mathcal{D}= \partial \Omega \setminus {\cal N} $ and $F$ be satisfying assumptions  (H1)--(H4).   Assume that $\beta:\R\to\R$ 
is a non decreasing continuous  function,  such that $\beta(0)=0$, and  let $f$ and $g$ be continuous in $\overline{\Omega}$, with $f\leq g$.  
If $u\in USC(\overline \Omega)$ satisfies 
                       $$\left\{\begin{array}{c}
                       -|\nabla u |^\alpha F( x, \nabla u, D^2 u) + \beta (u) \leq f \  \hbox{in} \ \Omega\\[1ex]
                       u\leq 0\ \hbox{ on} \ \mathcal{D}, \  \frac{\partial u}{ \partial \vec n } \leq 0\ \hbox{ on} \ {\cal N}
                       \end{array}\right.$$
and $v\in LSC(\overline \Omega)$  satisfies 
                        $$\left\{\begin{array}{c}
                       -|\nabla v |^\alpha F( x, \nabla v, D^2 v) + \beta (v) \geq g\  \hbox{ in} \ \Omega \\[1ex]
                       v\geq 0\ \hbox{on} \ \mathcal{D},\   \frac{\partial v}{\partial \vec n } \geq 0\ \hbox{on} \ {\cal N}
                       \end{array}\right.$$
and if either $\beta$ is strictly  increasing   or   $f< g$ in $\overline{ \Omega}$, then
                       $$u \leq v \ \hbox{ in} \ \overline{\Omega}\, .$$
                       \end{theo}
                       
                        \begin{proof} 
 The  proof closely  follows  the proof of Theorem 4.9 in \cite{Pa}, so we just present a sketch of it.
 
  Let us suppose by contradiction that $\max_{\overline \Omega} (u-v) >0$ and let us                       
introduce,  for  $q > \sup \left\{ 2, \frac{\alpha+2}{\alpha+1}\right\}$ and $L= \frac{q}{R}$, where $R$ is the uniform radius of  exterior spheres to $\Omega$,  the function
$$ \psi_j(x, y) =    u(x)-v(y)-\frac{j}{q} |x-y|^q e^{-L(d(x)+ d(y))}:=  u(x)-v(y)-\varphi(x,y)\, ,$$
with $j>0$ parameter devoted to tend to infinity. 

Since
$$
\max_{\overline{\Omega}^2} \psi_j(x,y)\geq \max_{\overline{\Omega}} (u(x)-v(x))\, ,
$$
 the maximum of $\psi_j$ on $\overline{\Omega}^2$ is also positive, and,  if $(x_j, y_j)$ denotes a maximum point for $\psi_j$,  then $j |x_j-y_j|^q \rightarrow 0$. Moreover,  there exists a maximum point $\bar x\in \overline{\Omega}$ for $u-v$ such that, up to a subsequence, 
 $$(x_j, y_j) \rightarrow (\bar x, \bar x)\, ,\quad  \psi_j(x_j,y_j)\to \max_{\overline \Omega}(u-v)\, ,\quad 
u(x_j) \rightarrow u( \bar x), \ v( y_j) \rightarrow v( \bar x).$$
Since $\max_{\overline{\Omega}}(u-v)=u(\bar x)-v(\bar x)>0$, it follows that $\bar x \notin  \mathcal{D}$. Let us select $r>0$ such that $B( \bar x, r) \cap \partial \Omega\subset {\cal N} $ if $\bar x \in {\cal N}$, and $B( \bar x, r) \subset \Omega$  if $\bar x \in \Omega$. For $j$ large enough, we have that $(x_j, y_j) \in B(\bar x, r)^2$. 

By Theorem \ref{lipN} and Remark \ref{analogies}, one has 
                          $$\frac{j}{q} |x_j-y_j |^qe^{-L(d(x_j)+ d(y_j))} + \max_{\overline{B(\bar x,r)\cap \Omega}} ( u-v) \leq u(x_j)-v( y_j) \leq  
\max_{\overline{B(\bar x,r)\cap \Omega}} ( u-v) + c  | x_j-y_j|$$
  which implies that 
\begin{equation}\label{gradq}
 j |x_j-y_j |^{q-1}\leq C\, . 
 \end{equation}
 Computing the derivatives of $\varphi$  one gets 
  $$\begin{array}{c}
   \nabla_x \varphi (x,y) = (j |x-y|^{q-2} (x-y) -\frac{Lj}{q}  |x-y|^q \nabla d(x))e^{-L(d(x)+d(y))}\, ,\\[2ex]
\nabla_y \varphi (x,y) = (-j |x-y|^{q-2}(x-y) -\frac{Lj }{ q} |x-y|^q \nabla d(y))e^{-L(d(x)+ d(y))}\, .
\end{array}$$
 Denote $p_j := \nabla_x \varphi (x_j,y_j), r_j = \nabla_y \varphi (x_j,y_j)$.
                              Since $x_j-y_j\rightarrow 0$, one has 
                              $|p_j| , |r_j| \geq C j |x_j-y_j| ^{q-1}$, for $j$ large enough, and 
                               arguing as in \cite{Pa}, one can suppose that $x_j\neq y_j$. Moreover, arguing as in the proof of Proposition \ref{holderN}, we get that both $x_j$ and $y_j$ belongs to $\Omega$.

By applying as usual Ishii's Lemma,  \cite{CIL}, we obtain that, for all $\epsilon >0$, there exist $X_j, Y_j\in \mathcal{S}$ such that $(p_j, X_j) \in \overline{J}^{2,+}u(x_j)$, $(-r_j, -Y_j) \in \overline{J}^{2,-} v(y_j)$ and 
$$-\left(\frac{1}{ \epsilon}+ |D^2 \varphi(x_j,y_j)|\right) I_{2N} \leq \left( \begin{array}{cc}
 X_j&0\\0& Y_j\end{array}\right) \leq  D^2 \varphi(x_j,y_j) + \epsilon (D^2 \varphi(x_j,y_j))^2.$$
Hence, we obtain
$$ |r_j|^\alpha  F( y_j, -r_j, -Y_j)-|p_j|^\alpha  F(x_j, p_j, X_j) +\beta (u(x_j))-\beta (v(y_j)) \leq f(x_j)-g(y_j)\, .$$                             
From here, following exactly the same arguments of \cite{Pa} and using \eqref{gradq}, we deduce, in the limit as $j\to \infty$,
$$  \beta (u(\bar x))-\beta ( v( \bar x))\leq f(\bar x) -g( \bar x)\, ,$$
which yields  a contradiction in both of the cases $f< g$ and $\beta$ is non decreasing , or $f\leq g$ and $\beta$ is strictly increasing.  \end{proof}

                 \subsection{An existence result}
In this section we provide the construction of  sub- and super-solutions under the 
 stronger assumption that $\Omega$ is of class ${ \cal C}^2$, and $\mathcal{D} = \overline{{\rm int}(\mathcal{D})}$, where   closure and   interior   are meant, as usual,  in the relative topology of $\partial \Omega$.
 
We consider, as before,  the problem
 $$\left\{\begin{array}{lc}
  -| \nabla u |^\alpha F(x, Du,  D^2 u) + \beta (u)= f& {\rm in} \ \Omega\\
 u= 0 \  {\rm on}\  \mathcal{D}, \frac{\partial u}{\partial \vec n } = 0\  {\rm on } \ {\cal N}&
\end{array}\right.$$
 where the datum $f$ is bounded and continuous in $\Omega$ and the zero order term $\beta:\R\to\R$ is  now assumed to satisfy, besides being continuous and  strictly increasing,
 \begin{equation}\label{beta}
 \lim_{t\to \pm \infty} \beta (t)=\pm \infty\, .
 \end{equation}
In order to construct a supersolution of the problem in $\Omega$, let us start by   constructing a super-solution in  $\Omega \cap B(z,r_z)$, where $z \in {\rm int}(\mathcal{D})$ and $r_z>0$ is chosen such that $r_z< \min \{ 1, {\rm dist}(z,\mathcal{N}), \delta \}$, with $\delta>0$ so  small   that   the distance from the boundary $d$ is ${ \cal C}^2$ in the set $\{ d<\delta\}$.

 Let $M>1$ be a large enough constant such that $\beta ( M) \geq  \|f\|_\infty$. 
 Let us set
    \begin{equation}\label{Gamma}
     \Gamma =  8 N \frac{A}{a} \left(  \frac{M}{r_z}  + \| D^2 d \|_\infty  \right)+ \left( \frac{2^{1+|\alpha|}}{a}  \|f\|_\infty  \right)^{\frac{1 }{2+ \alpha}}  ,
     \end{equation} 
      where, we recall,   $N$ is the dimension of the space, $a\leq  A$ are the ellipticity constants of the operator $F$, and 
$$
\| D^2d\|_\infty = \sup \{ |\lambda (x)|\, : \lambda (x) \hbox{ is an eigenvalue of }D^2d(x) \hbox{ and } d(x)<\delta \}\, .
$$

Let us further select a large constant $\kappa>1$ such that
                     \begin{equation} \label{kappa}
                     \beta  (\log ( 1+ \kappa) ) \geq N A (2\Gamma)^{\alpha^++2}\, ,
                     \end{equation}
 and let us finally set                    
 \begin{equation}\label{C}
  C =  (1+ \kappa) \Gamma\, .
 \end{equation}
 We consider the function        
 $$ \varphi_z(x) = \left\{ \begin{array}{lc}
                     \log ( 1+ Cd(x)) &  \hbox{ if } x\in  \Omega \cap B_{r_z} (z)\cap \{ Cd < \kappa\}\\
                      \log ( 1+ \kappa) &  \hbox{ if } x\in \Omega \cap B_{r_z} (z)\cap \{ Cd \geq  \kappa\}
                      \end{array}\right.$$ 
and we claim that
 $$   v_z(x)  = \varphi_z(x) + M \frac{ |x-z|^2 }{ r_z^2}$$
 is a super-solution in  $\Omega \cap B_{r_z} (z)$ in the sense of Definition \ref{defvis}. 
 Indeed,  if $Cd(x)<\kappa$, by direct computation we get
$$
\begin{array}{c}
\displaystyle \nabla v_z(x)= \frac{ C \nabla d(x)}{ 1+ Cd(x)}+ \frac{2M}{r_z^2}  ( x-z)\, ,\\[3ex]
\displaystyle D^2 v_z(x)= -C^2\frac{ \nabla d(x) \otimes \nabla d(x)}{(1+ Cd(x))^2}+ \frac{C D^2d(x)}{1+ Cd(x)} + \frac{2M}{ r_z^2} I_N\, .
\end{array}
$$
 Thus, by  using \eqref{Gamma} and \eqref{C} and observing that $\Gamma=\frac{C}{1+\kappa } <\frac{C}{1+C d(x)}$, we obtain
 $$ 
 | \nabla v_z(x)| \leq  \frac{C}{1+C d(x)} + 2 \frac{M}{r_z}\leq  \frac{C}{1+C d(x)}+ \Gamma \leq 2 \frac{C}{1+C d(x)}\, ,
 $$
 as well as
 $$  | \nabla v_z(x)|\geq \frac{C}{1+C d(x)}- \frac{1}{2} \Gamma \geq \frac{1}{2} \frac{C}{1+C d(x)}\, .
 $$
Hence
$$ | \nabla v_z(x) |^\alpha \geq 2^{-| \alpha|} \left(\frac{C}{1+C d(x)}\right)^{\alpha}\, .$$ 
Moreover,  by \eqref{Gamma}, it follows that
$$
\begin{array}{rl}
\displaystyle  \mathcal{M}^+_{a,A}( D^2  v_z(x))  \leq   & \displaystyle N\, A \left( \frac{C}{1+Cd(x)}\| D^2d\|_\infty +\frac{2M}{r_z^2}\right) -a\, \left( \frac{C}{1+Cd(x)}\right)^2\\[3ex]
\leq  & \displaystyle -\frac{a}{2}\left( \frac{C}{1+Cd(x)}\right)^2\, ,
\end{array}
 $$
so that
$$ | \nabla v_z(x) |^\alpha \mathcal{M}^+_{a,A}( D^2 v_z(x)) \leq -a 2^{-(1+| \alpha|)} \left( \frac{C}{1+ Cd(x)}\right)^{2+\alpha}\leq - a 2^{-(1+| \alpha|)} \Gamma^{2+\alpha}\, .$$
Since $\beta(v_z(x))\geq 0$, from \eqref{estiPucci} and  \eqref{Gamma} again,  we deduce
 $$- | \nabla v_z(x)|^\alpha F(x, \nabla v_z(x), D^2 v_z(x)) +\beta  ( v_z(x)) \geq  \|f\|_\infty\geq f(x)\, .$$
If $x$ belongs to the  set $\{ C d(x)=\kappa\}$, then neither there is a  $C^2$ function $\varphi$ such that $v_z-\varphi$ has a local minimum at $x$ nor $v_z$ is locally constant   around $x$, so that the  super-solution condition is certainly satisfied at $x$.
                             
\noindent Finally, if $C d(x)>\kappa$, we have on the  one hand  
$$|\nabla v_z(x)| = \left| \frac{2M}{r_z^2} (x-z)\right|\leq \frac{2M}{r_z}\, ,
$$
and, on the other hand,
$$|\nabla v_z(x)| \geq \frac{2M}{r_z^2} d(x)>\frac{2M}{r_z} \frac{\kappa }{C}\, .
$$
By using also that $\kappa >1$, we then obtain
$$
|\nabla v_z(x)|^{\alpha} \leq \left( \frac{2M}{r_z}\right)^\alpha \left( \frac{C}{\kappa}\right)^{\alpha^-}\leq \left( \frac{2M}{r_z}\right)^\alpha \left( 2 \Gamma\right)^{\alpha^-}\, .
$$
We further notice that
$$ \mathcal{M}^+_{a,A}( D^2v_z(x))\leq \frac{2NAM}{r_z^2}< 2NA \left( \frac{M}{r_z}\right)^2\, ,
$$
hence
$$|\nabla v_z(x)|^{\alpha} \mathcal{M}^+_{a,A}( D^2v_z(x))\leq 2^{1+\alpha^+} NA \left( \frac{M}{r_z}\right)^{2+\alpha} \Gamma^{\alpha^-}< \frac{NA}{2}(2\Gamma)^{2+\alpha^+}\, .
$$
 By using also \eqref{kappa}, we then deduce
 $$ -| \nabla v_z |^\alpha F(x,\nabla v_z, D^2 v_z)+ \beta ( v_z) > -\frac{NA}{2}(2\Gamma)^{2+\alpha^+} +\beta (\ln (1+\kappa)) > \|f\|_\infty\, ,$$
and the claim is proved.                                  
                  
 Next, we observe that the constant $M$ is also a super-solution in the whole domain $\Omega$, by the choice of $M$.  As a consequence,  $T_M(v_z)=\inf \{v_z, M\}$ is still a continuous super-solution in $\Omega\cap B_{r_z}(z)$, being  the infimum of two super-solutions. Moreover, $T_M(v_z)=M$ in a neighborhood of $\partial B_{r_z}(z)\cap \Omega$ and, by extending it  with the value  $M$, we obtain a continuous super-solution in the whole $\Omega$.  
 
  We then define 
                             $$ v (x)= \inf_{z\in \hbox{int}(\mathcal{D})} T_M(v_z)(x)\, .$$
As  it is well known, see \cite{CIL}, the lower semicontinuos envelope $v_*$  of $v$ is a lower semicontinuous super-solution. We observe that $v_*$ is not continuous on $\overline{\Omega}$, since 
 $v_*(x) = M$ for $x\in \mathcal{N}$ and 
$v_*(x) = 0$ on $\mathcal{D} = \overline{{\rm int} \mathcal{D}}$.  
                    
Symmetric arguments show that $(-v)^*=-v_*$ is a sub-solution of the considered problem.

We are now in the position of proving the following existence theorem
 \begin{theo}\label{existence} Given $f$ continuous and $\beta$ strictly increasing satisfying \eqref{beta}, 
 there exists a unique $u\in C(\overline \Omega)$ solution of
 $$\left\{\begin{array}{lc}
  -| \nabla u |^\alpha F(x, Du,  D^2 u) + \beta (u)= f& {\rm in} \ \Omega\\
 u= 0 \  {\rm on}\  \mathcal{D}, \frac{\partial u}{\partial \vec n } = 0\  {\rm on } \ {\cal N}&
\end{array}\right.$$
\end{theo}
Indeed we have now all the ingredients to use Perron's method adapted to the context, see \cite{Ishii1, BD1}, which yields  that the function
$$
u(x)=\sup \{ w(x)\, : \ w\in USC(\overline{\Omega} ) \hbox{ is a subsolution  and} -v_*\leq w\leq v_* \hbox{ in } \overline{\Omega}\}
$$
is the unique  continuous solution of the mixed boundary value problem. Let us point out that the standard argument in Perron's method for viscosity solutions can be applied due to the comparison result of Theorem \ref{compprinc}, and it gives that $u$ satisfies the Neumann boundary condition in the generalized viscosity sense. By Remark \ref{weakstrong}, the strong formulation is then recovered a posteriori. 
 
%%%%%%%%%%%%%%%%%%%%%%%%%%%%%%%%%%%%%%%%%%%%%%%%%%%%%%%%%%%%%%%%%%%%%%%%%%%%%%%%%%%%%%%%%%%%%%%%%%%%%%%%%%%%%%%%%%%%%%%%%%                       
  \section{Global H\"older regularity. }\label{section4}
In this section we make the  additional assumption that $\mathcal{D}\cap \overline{\cal N}$ is the union of $(N-2)$-dimensional manifolds of class ${ \cal C}^2$;
then, the following condition holds true:
\begin{itemize}
\item[(H5)]
For every $\bar x \in \mathcal{D}\cap \overline{\cal N}$, there exist  $r, \ r_0>0$,  and a $\mathcal{C}^2$-diffeomorphism $H_{\bar x, r}:B(\bar x,r)\to B(0,r_0)$ such that $H_{\bar x, r} (\bar x) = 0$, $H_{\bar x, r} ( B( \bar  x, r) \cap \Omega) = B(0,r_0) \cap \R^{N-1} \times \R^+$, and 
$$ \begin{array}{ll}
            H_{\bar x, r} ( \mathcal{D}\cap B(\bar x, r)) & \!\!\! = \{ x_{N-1} \leq  0\} \cap \{ x_N = 0\} \cap B(0,r_0),\\[1ex]
            H_{\bar x, r} ( {\cal N}  \cap B(\bar x, r)) & \!\!\! =  \{ x_{N-1} >0\} \cap \{ x_N = 0\} \cap B(0,r_0)\, .
\end{array}
$$
\end{itemize}
When no ambiguity arises,  we drop the indexes $\bar x, r$ for simplicity. We also say that (H5) holds uniformly if $r,\ r_0$ do not depend on $\bar x$ and, moreover, there exist universal constants $c, m>0$ depending only on $\Omega$ such that, in $B(\bar x,r)$, one has
$$
|\nabla H_{\bar x, r}| +|D^2 H_{\bar x, r}|+|\nabla H_{\bar x, r}^{-1}\circ H_{\bar x, r}|+|D^2 H_{\bar x, r}^{-1}\circ 
H_{\bar x, r}|\leq c
$$
and
$$|{\rm det} (\nabla H_{\bar x, r})|\geq m\, .
$$
We now proceed with the construction of some  barrier functions which will be used to get H\"older type estimates near points of $\mathcal{D}\cap \overline{\mathcal{N}}$.

 Let us start with a result in dimension $N = 2$, where, for  points $(x,y)\in \R^2$  we use the polar coordinates $r(x,y)=\sqrt{x^2+y^2}$ and    
 $$\theta (x,y)= \left\{ \begin{array}{lc}
       \arccos\frac{x}{ r} & \hbox{ if} \ y\geq 0\\[1ex]
        2\pi-\arccos\frac{x}{ r}  & \hbox{ if} \ y < 0
        \end{array}\right.$$
 \begin{lemme}\label{super2}
For $r_0 >0$ and $ 0<   \theta_1<   2\pi$ let us consider   the  conical domain
         $$ \mathcal{C}( r_0, \theta_1)  :=  \{r<  r_0, \ \theta \in (0, \theta_1)\}\, ,$$
and let further $\nu =(\nu_1,\nu_2)$ be a  vector field defined on $\{ \theta=0\}$ such  that $\nu_2 <0$.  
 
 Then, there exist  $\gamma \in ]0,1[$  and   a ${ \cal C}^2$ function  $\varphi:[0,  \theta_1]\to \R$, satisfying $\varphi \geq c>0$ for a universal constant $c>0$,  such that $w = r^\gamma  \varphi ( \theta)$ satisfies
    $$\left\{\begin{array}{l}
   -|\nabla w|^\alpha   F(x,  Dw, D^2 w) \geq c\, r^{ \gamma-2+ ( \gamma-1) \alpha} \ \hbox{ in  } \ \mathcal{C}( r_0, \theta_1),  \\[1ex]
   \displaystyle  \frac{\partial w}{\partial \nu} \geq  -c\,  \nu_2  r^{ \gamma-1}  \  \hbox{ on}\  \theta = 0\, .\\
     \end{array} \right.$$
           \end{lemme}
\begin{proof}
We denote 
       $X= ( x,y),  X^\perp = ( -y, x)$. Then, we observe that
       $$ \nabla \theta =\frac{ X^\perp }{r^2}, \quad        D^2 \theta = -\frac{X \otimes X^\perp+ X^\perp \otimes X}{r^4}.$$
Hence,   $D^2 \theta$ has eigenvalues $\pm \frac{1 }{ r^2}$. 
By a direct computation, one has  
         $$\nabla w =r^\gamma  \varphi^\prime \nabla \theta + \gamma r^{ \gamma-2} \varphi X, $$
          and, since $X$ and $\nabla \theta$ are orthogonal, 
          $|\nabla w | = r^{ \gamma-1} ((\varphi^\prime)^2 + \gamma ^2\varphi ^2)^{\frac{1}{2}}$. 
Moreover,
$$    
\begin{array}{rl}
          D^2 w = &\!\!\! r^\gamma  \varphi^{ \prime \prime} \nabla \theta \otimes \nabla \theta + \varphi^\prime ( \theta) \left( r^\gamma D^2 \theta + \gamma r^{ \gamma-2} ( X \otimes \nabla \theta + \nabla \theta \otimes X)\right)\\[2ex]
          &\!\!\! + \, \gamma r^{\gamma -2} \varphi \left( I+ ( \gamma-2) \frac{X \otimes X}{r^2}\right) \\[2ex]
          = &\!\!\! a r^{\gamma -4} \left\{ \varphi^{ \prime \prime}X^\perp \otimes X^\perp+(\gamma-1)\varphi^\prime (X^\perp \otimes X+X\otimes X^\perp) +\gamma\varphi(r^2I+(\gamma-2)X \otimes X )\right\}
\end{array}
$$
By conditions (H1) and (H2), and by assuming $\varphi>0, \varphi', \varphi''\leq 0$, it follows that
$$\begin{array}{ll}
F (x,\nabla w, D^2 w)&\!\!\! \leq \mathcal{M}^+_{a,A} (D^2w)\\[1ex]
           & \!\!\!  \leq   r^{\gamma-2}\left( a \varphi^{ \prime \prime} - (A-a) (1-\gamma) \varphi^\prime +  \gamma ( A +a( \gamma-1))\varphi \right)\, . 
            \end{array}
$$
We choose $\varphi$ satisfying
$$ \frac{a }{2} \varphi^{\prime \prime} -(A-a) (1-\gamma) \varphi^\prime +  \gamma ( A+a( \gamma-1))\varphi =0\, ,$$
  that is of the form
  $$\varphi =C_1 e^{\sigma_1 \theta} + C_2 e^{\sigma_2\theta}$$
  with
  $$
\sigma_{1,2} = \frac{ (A-a) (1-\gamma) \pm \sqrt{( A-a)^2(1-\gamma)^2-2a \gamma ( A+ (\gamma-1) a)}}{a}\, .$$
Note that, for $\gamma$ small enough, $\sigma_1 > 0$ and it is  close to $2 \frac{A-a}{a} (1-\gamma)$, and $\sigma_2  = c\,  \gamma  + O ( \gamma^2)$ and  close to zero.  Let us fix some $\kappa>2$, assume that  $\gamma$ is small enough in order that 
               $\kappa \frac{\sigma_2}{\sigma_1} e^{8\frac{(A-a)}{a} \pi} \leq \frac{1}{ 2}$. 
           Then, we choose $C_1 = -\kappa \frac{\sigma_2 }{\sigma_1}$ and $C_2 = 1$. Hence, we obtain
$$\varphi (\theta)= -\kappa \frac{\sigma_2 }{\sigma_1}e^{\sigma_1\theta} + e^{ \sigma_2 \theta}  \geq -\frac{ \kappa \sigma_2 }{\sigma_1} e^{8\frac{(A-a)}{a} \pi} +1\geq \frac{1}{2}\, ,$$
whereas 
              $$ \varphi^\prime(\theta) = -\kappa \sigma_2  e^{\sigma_1 \theta} + \sigma_2 e^{ \sigma_2 \theta} \leq -(\kappa-1) \sigma_2   e^{ \sigma_2 \theta}<0\, .$$
Furthermore, one has                             
               $$\varphi^{\prime \prime} = -\kappa \sigma_2 \sigma_1 e^{ \sigma_1 \theta} + \sigma_2^2 e^{ \sigma_2 \theta} \leq \sigma_2 ( \sigma_2-\kappa \sigma_1 ) e^{ \sigma_2 \theta}  < -(\kappa-1)\sigma_1\sigma_2  e^{ \sigma_2 \theta}<0\, . $$
Hence, we deduce                 
            $$ F(x,  \nabla w, D^2 w) \leq \frac{a}{2} r^{\gamma-2}\varphi^{\prime \prime} \leq -\frac{a \sigma_1\sigma_2}{2}  r^{\gamma-2}\, , $$
and, since
$$ c_1 r^{ \gamma-1}\leq |\nabla w|=r^{ \gamma-1} \sqrt{(\varphi^\prime)^2 + \gamma ^2\varphi ^2}\leq c_2 r^{ \gamma-1}$$
for positive universal constants $c_1$ and $c_2$, in conclusion we  get
           $$| \nabla w|^\alpha   F(x,  Dw, D^2 w) \leq -c r^{ \gamma-2+ ( \gamma-1) \alpha}. $$
Moreover, for any given $\nu=(\nu_1,\nu_2) \in \R^2$ with $\nu_2<0$ and for $\theta=0$, one has           
 \begin{eqnarray*}
                 \nabla w  \cdot \nu  &=& r^{ \gamma-1} (\gamma \varphi (0)\nu_1 +  \varphi^\prime(0) \nu_2) \\
                 &=&  r^{ \gamma-1} \left(\gamma \nu_1 ( -\kappa \frac{ \sigma_2 }{\sigma_1} + 1) - \nu_2 (\kappa -1)  \sigma_2 \right)
                 \geq -c \nu_2\gamma   r^{ \gamma-1}
                 \end{eqnarray*}
  for $\gamma$ sufficiently small. \end{proof}

 In  the next result we extend the above construction from dimension $2$ to dimension $N$, where for points $x=(x',x_{N-1},x_N)\in \R^N$ we use the cylindrical coordinates $x'\in \R^{N-2}$, $r=\sqrt{x_{N-1}^2+x_N^2}$, $\theta=\theta(x_{N-1},x_N)$.
\begin{theo}\label{superso} For $r_0>0$ and $\theta_1\in (0,2\pi )$, let
$$\tilde {\mathcal{C}}(r_0,\theta_1)= \{ |x'|< r_0\} \times \mathcal{C}(r_0,\theta_1)\, ,$$
and let $\nu=(\nu_1,\ldots,\nu_N)$ be a vector field defined on $\{ \theta=0\}$ satisfying $\nu_N \leq -m<0$. There exist
$\gamma\in ]0,1[$,    $C,\ r_0>0$,  depending only on the dimension and the ellipticity constants, and a function $\tilde w:\overline{\tilde{\mathcal{C}}(r_0,\theta_1)}\to \R$, satisfying 
           $$\frac{1}{ C} (r^\gamma+ |x^\prime |^2) \leq \tilde{w}\leq C(r^\gamma + | x^\prime |^2)\, ,$$   
such that
       $$\left\{\begin{array}{l}
 \displaystyle   -|\nabla \tilde{w}|^\alpha   F(x,  \nabla \tilde{w}, D^2 \tilde{w}) \geq\frac{1}{ C}\,  r^{ \gamma-2+ ( \gamma-1) \alpha} \ \hbox{ in  } \   \tilde{{\mathcal C}}( r_0, \theta_1) \\[1ex]
 \displaystyle     \frac{\partial \tilde{w}}{\partial \nu} \geq \frac{1}{ C}\,   r^{ \gamma-1} \ \hbox { on } \ \theta = 0  
     \end{array} \right.$$
\end{theo}
 \begin{proof}
            It is sufficient to consider in the plane $(x_{N-1}, x_N)$ the function $w$ constructed in Theorem \ref{super2}, relabelling the coordinates. Then, let us set $\tilde{w}= w+ C |x^\prime |^2$. One has
$$|\nabla \tilde{w} |^2 \leq  C(|x^\prime|^2 + r^{2(\gamma-1)} )\leq    C |x|^{2(\gamma-1) }\, ,$$
as well as
               $$| \nabla  \tilde{w}|^2 \geq C r^{2(\gamma-1)} \geq   C  |x|^{2(\gamma-1) }\, .$$
               Furthermore, 
               $$D^2  \tilde{w} = \left( \begin{array}{cc} 
               2 C\, I_{N-2} & 0\\
                0 &  D^2w\end{array} \right)\, .$$
An easy computation shows that $\tilde{w}$ satisfies the statement.
                                The details are left to the reader.                \end{proof}

\begin{rema}
{\rm  Under the current assumption (H5), we will use Theorem \ref{superso} with $\theta_1 = \pi$.}
                                            \end{rema}
                              
Next, we analyze the  change of equation under the diffeomorphism $H$  appearing in condition (H5). In the sequel, for 
$y\in B(0,r_0) \cap \R^{N-1} \times \R^+$, $p\in \R^N$ and $M\in \mathcal{S}$, we set        
$$F_H( y,p, M) = F\left(G(y), p \nabla H (G(y)), ^t \nabla H(G(y))  : M: \nabla H(G(y))+ p\, D^2H(G(y))\right).$$ 
                                         
 \begin{lemme}\label{change}
Assume condition (H5) holds true, and let $H=H_{\bar x,r}$, $G=H^{-1}$. If $u\in USC (\overline{\Omega \cap B(\bar x, r)})$ satisfies
$$\left\{\begin{array}{l}
 -  |\nabla u |^\alpha F(x, \nabla u,  D^2 u) \leq f \ \hbox{  in } \ \Omega\cap B( \bar x, r),\\[1ex]
 \displaystyle u\leq  0 \  \hbox{ on } \ \mathcal{D} \cap B( \bar x, r),\quad  \frac{\partial  u}{\partial \vec n } \leq 0 \ \hbox{ on } \ {\cal N}  \cap B( \bar x, r),
 \end{array}\right.$$
 then, the function $v(y)=u(G(y))$ satisfies
$$\left\{\begin{array}{l}
- | \nabla v (\nabla H (G(y))|^\alpha F_H( y,\nabla v, D^2 v) \leq f(G(y)) \ \hbox{ in }  B(0, r_0)\cap \R^{N-1} \times \R^+  ,\\[1ex]
\displaystyle \  v\leq  0 \  \hbox{ on } \ H(\mathcal{D}\cap B( \bar x, r))\, ,\quad
\frac{\partial v}{\partial \nu}  \leq 0 \ \hbox{ on } \ H({\cal N}\cap B( \bar x, r))\, ,
  \end{array}\right.$$
 for a smooth vector field $\nu$ defined on $H({\cal N}\cap B( \bar x, r))$ and satisfying $\nu_N<0$.
                                               
An analogous  result holds for super-solutions.
 \end{lemme}
 \begin{proof} If $u$ is of class $\mathcal{C}^2$, then, using the convention of summation over repeated indexes, from the identity $u(x)=v(H(x))$ we deduce, for $1\leq i,j\leq N$,
 $$
 \frac{\partial u}{\partial x_i}(x)=  \frac{\partial v}{\partial y_k}(H(x))\, \frac{\partial H_k}{\partial x_i}(x)\, ,
 $$
 as well as
 $$
 \frac{\partial^2u}{\partial x_i\partial x_j}(x)= \frac{\partial H_l}{\partial x_j}(x) \frac{\partial^2v}{\partial y_k\partial y_l}(H(x)) \frac{\partial H_k}{\partial x_i}(x) +  \frac{\partial v}{\partial y_k}(H(x)) \frac{\partial^2H_k}{\partial x_i\partial x_j}(x).
 $$
 Moreover,  by assumption (H5), the last component $H_N$ of the map $H$ satisfies $H_N(x)>0$ for $x\in \Omega\cap B(\bar x,r)$ and $H_N(x)=0$ for $x\in \partial \Omega\cap 
 B(\bar x,r)$. Hence, one has
  $$
  \vec n(x)=-\frac{\nabla H_N(x)}{|\nabla H_N(x)|}\quad \hbox{ for } x\in \partial \Omega \cap B(\bar x,r)\, ,
  $$
 so that, on $B(\bar x,r)\cap \mathcal{N}$, we have
  $$
  \nabla u(x)\cdot \vec n(x)= \nabla v (H(x)) \cdot \nu (H(x))\, ,
  $$
  with
  $$
  \nu (H(x))= \nabla H(x) \cdot \left( -\frac{\nabla H_N(x)}{|\nabla H_N(x)|}\right)
  $$
  and
  $$
  \nu_N(H(x))= -|\nabla H_N(x)|<0\, .
  $$
The extension of the above result to the case when $u$ is merely a viscosity subsolution  immediately follows from the viscosity formulation.
 \end{proof}
 
 \begin{rema}
 {\rm If assumption (H5) is uniformly satisfied, then the vector field $\nu$ given in the above Lemma   satisfies
 $$
 \nu_N\leq -m<0\quad \hbox{ in } H(\mathcal{N}\cap B(\bar x,r))\, ,
 $$
 for a universal constant $m>0$ depending only on $\Omega$ and $\mathcal{N}$.}
 \end{rema}      
  
  We can now prove local H\"older type inequalities for sub- and supersolutions near points of $\mathcal{D}$.            
    \begin{theo}\label{localholder}
 Assume that condition {\rm (H5)}  is uniformly satisfied and   
let  $u\in USC(\overline \Omega)$  be  a  solution of 
$$\left\{
\begin{array}{l}
-| \nabla u |^\alpha F(x,\nabla u, D^2 u ) \leq  f\  \hbox{  in } \ \Omega,\\[1ex]
\displaystyle  u \leq  0\ \hbox{ on} \ \mathcal{D}, \quad \frac{\partial u}{\partial \vec n }  \leq  0 \  \hbox{ on } \ \mathcal{N}\, ,
 \end{array}\right.$$
                   with $f$ bounded  in $\Omega$. 

Then, there exists  $\gamma\in (0,1)$ such that, for any $r>0$ sufficiently small, there exists  $C_r>0$ satisfying, for all $\bar x \in \mathcal{D}$ and $x\in \overline{ \Omega} \cap B( \bar x, r)$, 
                    \begin{equation}\label{Cgamma}
                    u(x) \leq C_r |x-\bar x|^\gamma\, .
                    \end{equation} 
An analogous result holds true for supersolutions.
\end{theo}
\begin{proof} We give the proof when $\bar x\in \mathcal{D}\cap \overline{\mathcal{N}}$. For the case $\bar x\in \mathcal{D}
\setminus \overline{\mathcal{N}}$, see \cite{BDCocv}.
                               
We use the change  of variable given by the diffeomorphism $H$ and apply Lemma \ref{change}. 
It is sufficient to prove that the function $v(y)= u(G(y))$ satisfies in $B( 0, r_0) \cap  \R^{N-1} \times \R^+$
\begin{equation}\label{holderv}
v(y)  \leq C |y|^\gamma,
\end{equation}
 since, then, estimate \eqref{Cgamma} for $u$ automatically follows.   
 
In order to prove \eqref{holderv}  we just use the comparison principle.    Let $\tilde{w}$ be given by Theorem \ref{superso}, with $\theta_1=\pi$ and $r_0$ so small that, for $r < r_0$, one has $C r^{ \gamma-2+ ( \gamma-1) \alpha } > |f|_\infty$. Then, Theorem  \ref{compprinc}  applied  to  $v$ and   $\tilde{w}$ yields
                      $v\leq \tilde{w}$ in $B(0,r_0)\cap( \R^{N-1}\times \R^+)= \tilde{{\mathcal C}}( r_0, \pi)$, so that, again by Theorem \ref{superso},
                      $$v(y)\leq C\left( (y_{N-1}^2+y_N^2)^{\gamma/2} + | y^\prime |^2\right)\leq C |y|^\gamma\, .$$
                            \end{proof}
                            
                            \begin{rema}
 {\rm  Close to points  $\bar x\in {\rm int} (\mathcal{D})$  ( as well as  $\bar x \in {\cal N}$), one can prove  a better estimate, namely  a Lipschitz-type estimate, but with a constant depending on   ${\rm dist}( \bar x,\overline{\cal N})$, and  so not helpful  for   global H\"older  estimates.}
                            \end{rema} 
 
 We can now state and prove the main result of this section.                      
  \begin{theo}\label{globalholder}
Assume that condition {\rm (H5)} holds uniformly and let    $u\in USC(\overline{\Omega})$ and $v\in LSC(\overline{\Omega})$ be satisfying  $\max_{\overline{\Omega}} ( u-v) \geq 0$ and , respectively, 
 $$\left\{\begin{array}{l}
 -| \nabla u |^\alpha F(x, \nabla u, D^2 u) \leq f \  \hbox{ in } \ \Omega\, ,\\[1ex]
  \  u\leq  0 \  \hbox{ on } \ \mathcal{D}, \  \frac{\partial  u}{\partial \vec n } \leq 0 \  \hbox{ on } \ {\cal N}\, ,                                             \end{array}\right.$$  
 $$\left\{\begin{array}{l}
 - |\nabla v |^\alpha F(x, \nabla v, D^2 v) \geq g \  \hbox{ in } \ \Omega\, ,\\[1ex]
  \  v\geq  0 \  \hbox{ on } \ \mathcal{D}, \  \frac{\partial  v}{\partial \vec n} \geq 0 \  \hbox{ on } \ {\cal N}\, ,                                             \end{array}\right.$$  
with $f$ and $g$ continuous and bounded in $\Omega$. Then, for $\gamma\in (0,1)$ as in Theorem \ref{localholder} and   $M>0$  depending on the ellipticity constants of $F$, on $|f|_\infty$ and $|g|_\infty$, one has, for all $(x,y) \in \overline{ \Omega}^2$, 
 \begin{itemize}
 \item[(i)]  if $\max_{\overline{\Omega}} ( u-v) >0$, then
 $$ u(x)-v(y) \leq \max_{\overline{\Omega}} ( u-v) + M | x-y|\, ;$$
   \item[(ii)] if $\max_{\overline{\Omega}} ( u-v) = 0$, then 
  $$ u(x)-v(y) \leq  M | x-y|^\gamma\, . $$
 \end{itemize}
   In particular,  the solution of  the mixed boundary value problem \eqref{defmix} 
  is $\gamma$-H\"older continuous in $\overline{\Omega}$, with constant $M$. 
   \end{theo}
 \begin{proof}
  Assume  first that $\max_{\overline{\Omega}} ( u-v) >0$, and let us define  
 $$\psi(x,y) = u(x)-v(y)-\max_{\overline{\Omega}}(u-v) - M_\delta | x-y| e^{-L((d(x)+ d(y))}\, ,$$ 
with $L=\frac{2}{R}$, where $R>0$ is the uniform radius of exterior tangent spheres to $\partial \Omega$,  and   $M_\delta  =  \frac{2 (\sup  u-\inf v) e^{2L{\rm diam}(\Omega)}}{\delta}$, where $\delta>0$ has to be fixed suitably small.

  By arguing, as usual, by contradiction, we   suppose that 
 $\max_{\overline{\Omega}^2} \psi(x, y) >0$, so that it is achieved at a point $(\bar x, \bar y)$ satisfying
  $|\bar x-\bar y | < \delta$. 
  
  Let us check that neither $\bar x$ nor $\bar y$   belongs to $\mathcal{D}$. Indeed,  if, for instance,  $\bar x \in \mathcal{D}$, by the lower semicontinuity of $v$ we  can take $\delta$ small enough in order that   
 $$ u(\bar x)-v( \bar y) \leq -v(\bar y)\leq \frac{1}{2}\max_{\overline{\Omega}} (u-v),$$
 and, therefore,
 $$
 \psi(\bar x, \bar y)\leq  -\frac{1}{2}\max_{\overline{\Omega}} (u-v) - M_\delta | \bar x-\bar y| e^{-L((d(\bar x)+ d(\bar y))}<0\, ,
 $$
 yielding a contradiction. Analogously, we obtain that $\bar y\notin \mathcal{D}$.    The rest of the proof runs  as the proof of  Theorem  \ref{lipN}, and, for $\delta>0$ small enough, we get the conclusion.
                                 
  Suppose now that $\max_{\overline{\Omega}} (u-v)= 0$.  In this case, we set
  $$\psi(x,y) = u(x)-v(y)- M _\delta | x-y|^\gamma e^{-L((d(x)+ d(y))}$$
 and we select $\delta$ so small that $M_\delta  > 2 e^{L {\rm diam}(\Omega)} C_{2\delta}$,   where $C_{2\delta}$ is the constant  given by (\ref{Cgamma}).  We claim that, also in this case, neither $\bar x$ nor $\bar y$ belongs to $\mathcal{D}$. Indeed, assuming by contradiction that $\bar x \in  \mathcal{D}$, then, from  Theorem \ref{localholder}  it follows that
   $$ v( \bar y) \geq - C_{2\delta} |\bar  x
-\bar y |^\gamma\, , $$ 
 and, therefore,
  $$
   u(\bar x)-v( \bar y) \leq 0-v( \bar y)\\ 
     \leq  \frac{M_\delta }{2} e^{-L d(\bar y)} | \bar x-\bar y |^\gamma $$ 
 which again contradicts   $\psi( \bar x, \bar y)>0$. 
       An analogous argument treats  the case $\bar y \in{ D}$.  
        The rest of the proof follows the lines of the proof of Proposition  \ref{holderN}.                          
     \end{proof}
 The above theorem immediately yields  the following compactness result.   
  \begin{cor}\label{compact}
Assume that condition {\rm (H5)} holds uniformly and let  $(u_n)_n$ be a sequence of solutions of 
   $$ \left\{\begin{array}{l}
 -| \nabla u_n |^\alpha F( x,  \nabla u_n, D^2 u_n)  = f_n \ \hbox{ in } \ \Omega\, , \\
    u_n = 0 \ \hbox{ on} \ \mathcal{D},\ \frac{\partial  u_n}{\partial \vec n}  = 0\  \hbox{ on } \mathcal{N}\, . 
  \end{array}\right.$$
   Suppose that  $(u_n)_n$ is uniformly bounded, and that $f_n$ converges  locally uniformly in $\Omega$  to a continuous and bounded function $f$. Then,  up to a subsequence,  $(u_n)_n $   converges locally uniformly  in $\Omega$ to a  solution of 
 problem \eqref{defmix}.
  \end{cor} 

Finally,  we can also  easily derive the following existence result, which complements Theorem \ref{existence} for equations not having zero order terms. In this case, we need to add the positive homogeneity assumption
\begin{equation}\label{homo}
F(x, t\, p, t\, M)= t\, F(x,p, M)\qquad \hbox{ for all } x\in \Omega\, ,\ p\in \R^N\, ,\ M\in {\mathcal S}_N \hbox{ and } t\geq 0\, .
\end{equation}
\begin{theo}\label{existence1} Under assumptions {\rm (H5)} and \eqref{homo}, given $f$ continuous and bounded in $\Omega$, 
 there exists  $u\in C(\overline \Omega)$ solution of
 $$\left\{\begin{array}{lc}
  -| \nabla u |^\alpha F(x, Du,  D^2 u) = f& {\rm in} \ \Omega\, ,\\[1ex]
\displaystyle  u= 0 \  {\rm on}\  \mathcal{D}, \frac{\partial u}{\partial \vec n } = 0\  {\rm on } \ {\cal N}\, .&
\end{array}\right.$$
\end{theo}
\begin{proof} For every $n\geq 1$, by Theorem \ref{existence}, there exists a unique  $u_n\in C(\overline \Omega)$   solution to
$$\left\{\begin{array}{lc}
\displaystyle   -|\nabla u_n |^\alpha F(x, Du_n,  D^2 u_n) + \frac{|u_n|^\alpha}{n} u_n= f& {\rm in} \ \Omega\, ,\\[1ex]
\displaystyle   u_n= 0 \  {\rm on}\  \mathcal{D}, \ \frac{\partial u_n}{\partial \vec n } = 0\  {\rm on } \ {\cal N}\, .&
\end{array}\right.$$
If we show that the sequence $(u_n)_n$ is uniformly bounded, then, by Theorem \ref{globalholder}, we can deduce the local uniform convergence, up to a subsequence, to a function $u\in C(\overline \Omega)$. By the stability properties of viscosity solutions, see \cite{CIL}, it follows that $u$ is a solution, as claimed. 

In oder to prove the boundedness of $(u_n)_n$, we argue by contradiction, and we assume that, for a subsequence still denoted by $(u_n)_n$, one has $\| u_n\|_{\infty}\to \infty$. Then, by \eqref{homo}, the normalized function $\displaystyle v_n\, := \frac{u_n}{\|u_n\|_{\infty}}$ satisfies
$$\left\{\begin{array}{lc}
\displaystyle   -|\nabla v_n |^\alpha F(x, Dv_n,  D^2 v_n) + \frac{|v_n|^\alpha}{n} v_n= \frac{f}{\|u_n\|_{\infty}^{\alpha+1}}& {\rm in} \ \Omega\, ,\\[1ex]
\displaystyle   v_n= 0 \  {\rm on}\  \mathcal{D}, \ \frac{\partial v_n}{\partial \vec n } = 0\  {\rm on } \ {\cal N}\, .&
\end{array}\right.$$
Theorem \ref{globalholder} and the stability theory now apply to the bounded sequence $(v_n)_n$, yielding that, up to a subsequence, $(v_n)_n$ locally uniformly converges to a function $v\in C(\overline \Omega)$ satisfying
$$\left\{\begin{array}{lc}
\displaystyle   -|\nabla v |^\alpha F(x, Dv,  D^2 v) = 0& {\rm in} \ \Omega\, ,\\[1ex]
\displaystyle   v= 0 \  {\rm on}\  \mathcal{D}, \ \frac{\partial v}{\partial \vec n } = 0\  {\rm on } \ {\cal N}\, .&
\end{array}\right.$$
On the other hand, the strong maximum principle of \cite{BDCocv} jointly with Hopf boundary lemma yields that the above problem has only the trivial solution $v=0$. This clearly contradicts the fact that $\| v\|_{\infty}=1$ and completes the proof.
\end{proof}

%%%%%%%%%%%%%%%%%%%%%%%%%%%%%%%%%%%%%%%%%%%%%%%%%%%%%%%%%%%%%%%%%%%%%%%%%%%%%%%%%%%%%%%%%%%%%%%%%%%%%%%%%%%%%%%%%%%%%%%%%%%%%%                                      
                          
\section{Appendix}\label{appendix}         
 In this section we give the proof of Theorem \ref{lipN}. In order to do so we need to first prove a  local H\"older-type  estimate for sub- and super-solutions, 
 up to the boundary,  near ${\cal N}$.  
\begin{prop}\label{holderN}
    Let $\bar x\in \mathcal N$ and assume that $B(\bar x,\rho) \cap \partial \Omega \subset {\cal N}$. 
    Suppose that $u$ is  $USC (\overline{B(\bar x,\rho) \cap  \Omega})$ and satisfies 
     $$\left\{ \begin{array}{c}
      -| \nabla u |^\alpha F(x, \nabla u,  D^2 u)  \leq f \  \hbox{in} \  B(\bar x,\rho) \cap  \Omega \\[1ex]
    \displaystyle    \frac{\partial u}{\partial \vec n} \leq 0 \  \hbox{on} \  B(\bar x,\rho) \cap \partial \Omega 
             \end{array}\right.$$
as well as  $v$ is $LSC(\overline{B(\bar x,\rho) \cap  \Omega})$ and satisfies 
       $$\left\{ \begin{array}{c}
      -| \nabla v |^\alpha F(x, \nabla v,  D^2 v)  \geq g \ \hbox{in} \ B(\bar x,\rho) \cap  \Omega\\[2ex]
 \displaystyle \frac{\partial v}{\partial \vec n} \geq 0 \  \hbox{on} \  B(\bar x,\rho) \cap \partial \Omega 
      \end{array}\right.$$
 Suppose in addition that $f$ and $g$ are  bounded. Then   for any $\gamma<1$ and $r<\rho$ there exists $c_{r, \gamma}>0$, depending on $\Omega$, $N$, the ellipticity constants of $F$, $\|f\|_\infty,\ \|g\|_\infty$ and on $\sup u$, $\inf v$,  such that, 
       for all $(x, y)\in (B(\bar x,r) \cap \overline{ \Omega} )^2$, one has
       $$ u(x)-v( y) \leq \sup_{ \overline{ \Omega \cap B(\bar x,\rho)}}  ( u-v) + c_{r, \gamma} | x-y|^\gamma.$$
\end{prop}       
\begin{proof} 
Let us fix $r<\rho$, $x_0 \in B(\bar x, r)\cap \overline{\Omega}$, and let us introduce the function 
         $$\psi(x,y) = u(x)-v(y)-\sup_{ \overline{ \Omega\cap B(\bar x,\rho)}}  ( u-v) - (M | x-y|^\gamma + L_1 |x-x_0|^2+ L_1 |y-x_0|^2) e^{ -L(d(x)+ d(y))},$$
 for $M, L, L_1>0$ to be suitably chosen. We will prove that for $M, L, L_1$ large enough, independently  of $x_0$, but depending only on the ellipticity constants, on the supremum of
          $u-v$,  on $\|f\|_\infty,   \|g\|_\infty $, $\sup u$, $\inf v$ and finally  on $r$, the function $\psi$ satisfies
$$\psi(x,y) \leq 0\ \hbox{ in }    \overline{ B(\bar x,\rho)\cap  \Omega}\, .$$ 
         Once  this will have been proved,    taking $x= x_0$ and making $x_0$ vary  will provide  the result. 
         
Without loss of generality, we assume that $\rho$ is so small that the distance function $d$ is of class $C^2$ in $\overline{\Omega \cap B(\bar x, \rho)}$. Moreover,  
 if  $R>0$ is a uniform radius such that the   exterior sphere condition is satisfied in any compact portion of ${\mathcal N}$,  and  if $x\in {\mathcal N}$ and $y \in \overline{\Omega}$, then one has 
\begin{equation}\label{exter}
|x-y| ^2\geq 2R (x-y)\cdot  \nabla d(x). 
\end{equation}
         For the sequel we will   set
          $$ L =  \frac{2}{R}\, , $$
where $R>0$ is such that \eqref{exter} holds true for all $x\in {\mathcal N}\cap B(\bar x, \rho)$.

 Next,    we set 
$$\begin{array}{c}
\displaystyle S= \sup_{\overline{ \Omega \cap B(\bar x,\rho)}}u - \inf_{\overline{ \Omega \cap B(\bar x,\rho)}} v- \sup_{\overline{ \Omega \cap B(\bar x,\rho)}} (u-v)\, ,\\[3ex]
\displaystyle L_1 = \frac{4 S e^{2 L {\rm diam}(\Omega)}}{(\rho-r)^2}\, ,\\[2ex]
\displaystyle M =  \frac{S e^{2 L {\rm diam}(\Omega)}}{\delta^\gamma}\, ,
 \end{array} $$
 for $0<\delta<1$ to be chosen later small enough.
 
           Suppose now, by contradiction, that the supremum of $\psi$ on   $\overline{ \Omega \cap B(\bar x,\rho)}^2$ is strictly positive.  Then, it is achieved at some point $( \hat x, \hat y)$ belonging to  $( \Omega\cap B(\bar x,\rho) )^2$. Indeed, by the choice of $L_1$ and $M$, one has 
            $ (\hat x, \hat y) \in  B( x_0, \frac{\rho -r}{2})^2$, and $0<|\hat x-\hat y | < \delta$. In particular, both  $\hat x$ and $\hat y$ belong to  $B(\bar x, \frac{\rho +r}{2}) \cap \overline{ \Omega}$.  Let us further prove that neither $\hat x$ nor $\hat y$ belongs to $\partial \Omega$.
            
            \noindent  Indeed, we note  that, in particular, $\hat  x$ is a local maximum point of $u-\phi$, where
$$\phi (x)\, := \left( M|x-\hat y|^\gamma +L_1 |x-x_0|^2\right)e^{-L(d(x)+d(\hat y))}.$$
Assuming by contradiction that $\hat x\in \partial \Omega$, by a direct computation we get
$$
\begin{array}{rl}
\displaystyle \frac{\partial \phi}{\partial \vec n}(\hat x)  = e^{-L d(\hat y)} & \displaystyle \left[ M |\hat x-\hat y|^{\gamma -2} \left( L |\hat x-\hat y|^2-\gamma (\hat x-\hat y)\cdot \nabla d(\hat x)\right)\right.\\[2ex]
& \displaystyle \quad \left. +L_1\left(  L |\hat x-x_0|^2-2 (\hat x-x_0)\cdot \nabla d(\hat x)\right)\right]
\end{array}
$$
and therefore, by \eqref{exter} and the choice of $L$, we obtain
$$
\frac{\partial \phi}{\partial \vec n}(\hat x) \geq  e^{-L d(\hat y)}  \left[ M |\hat x-\hat y|^\gamma \left( L-\frac{\gamma}{2R}\right) +L_1 |\hat x-x_0|^2\left( L-\frac{1}{R}\right)\right]>0\, .
$$ 
Since $u$ satisfies $\frac{\partial u}{\partial \vec n } \leq 0$ on $\mathcal{N}$ in the viscosity sense, we    deduce that $\hat x\notin \mathcal{N}$. An analogous argument, applied to the super-solution $v$, shows that $\hat y\notin \mathcal{N}$.

Hence, both $\hat x$ and $\hat y$ belong to $\Omega \cap B(\bar x, \rho)$. Setting 
 $$\varphi(x,y) = (M | x-y|^\gamma + L_1 |x-x_0|^2+ L_1 |y-x_0|^2) e^{ -L((d(x)+ d(y))}\, ,$$
 Ishii's Lemma (see \cite{CIL}) yields that, for every $\epsilon >0$, there exist $X_\epsilon$ and $Y_\epsilon$ in $\mathcal{S}$ such that 
              $$( \nabla _x \varphi( \hat x, \hat  y), X_\epsilon) \in \overline{J^{2,+}} u( \hat  x), \ ( -\nabla_y \varphi(\hat  x, \hat  y), -Y_\epsilon) \in \overline{J^{2,-}} v( \hat  y)\, ,$$
with $X_\epsilon$ and $Y_\epsilon$ satisfying 
\begin{equation}\label{matrix}
  -\left( |D^2\varphi (\hat x,\hat y)|+\frac{1}{\epsilon}\right)  I_{2N}\leq 
    \left( \begin{array}{cc}         X_\epsilon& O \\
               O& Y_\epsilon\end{array} \right) \leq 
D^2 \varphi (\hat x, \hat y) +\epsilon \left( D^2 \varphi (\hat x, \hat y)\right)^2\, ,
\end{equation}
$I_{2N}$ denoting the identity matrix in $\R^{2N}$.

By a direct computation, one has
$$
\begin{array}{rl}
              \nabla_x \varphi(x,y) = e^{ -L(d(x)+ d(y))} & \!\!\! \left[  \gamma M | x-y |^{ \gamma-2} (x-y) + 2L_1 ( x-x_0)\right.\\[1ex]
             & -\left. L \nabla d(x) (M | x-y|^\gamma + L_1 |x-x_0|^2+ L_1 |y-x_0|^2) \right]
\end{array}
$$   
as well as
$$             
\begin{array}{rl}
\nabla_y \varphi(x,y) = e^{ -L(d(x)+ d(y))} & \!\!\! \left[ \gamma M | x-y |^{ \gamma-2} (y-x) + 2L_1 (y-x_0) \right.\\[1ex]
                 & -\left. L \nabla d(y) (M | x-y|^\gamma + L_1 |x-x_0|^2+ L_1 |y-x_0|^2) \right]\, .
\end{array}
$$                 
Hence,   for $\delta$ small enough depending on $r, \rho, R$ and $\gamma$, we have
\begin{equation}\label{Dx}
 \frac{ \gamma  M}{2} |\hat x-\hat y|^{ \gamma-1} e^{-2L {\rm diam } (\Omega)} \leq |\nabla_x \varphi (\hat x, \hat y)|,\  
|\nabla_y \varphi (\hat x, \hat y)|\leq \frac{3}{2} \gamma M  |\hat x-\hat y|^{ \gamma-1}.
\end{equation} 
Furthermore, we also have 
 \begin{equation}\label {Dx+Dy}
 |\nabla_x \varphi(\hat x, \hat y)+  \nabla_y \varphi(\hat x, \hat y)| \leq 2L_1 (\rho-r)+L\left(2M\delta^\gamma +L_1 (\rho-r)^2\right)=C\, .
                   \end{equation} 
 Here and in the following we denote by  $C$ a positive constant, which may vary from line to line, and which   depends on the data  of the problem and on $r$ and $\gamma$, but it is independent of  $\delta$.  
 
Estimates  \eqref{Dx}  in particular imply  that   $\nabla_x\varphi (\hat x, \hat y),\    \nabla_y\varphi (\hat x, \hat y)  \neq 0$. Since $u$ is a sub-solution and  $v$ a super-solution, we then get the two inequalities
 $$
 \begin{array}{c}
   -| \nabla_x \varphi(\hat x,\hat y) |^\alpha F(\hat x, \nabla_x \varphi(\hat x,\hat y),  X_\epsilon)  \leq f (\hat x)\, ,\\[2ex]
  -| \nabla_y \varphi(\hat x,\hat y) |^\alpha F(\hat y, -\nabla_y \varphi(\hat x,\hat y),  -Y_\epsilon)  \geq g (\hat y)\, .
  \end{array}
  $$
Hence, denoting simply by $(\nabla_x\varphi, \nabla_y \varphi)$ the gradient of $\varphi$ evaluated at $(\hat x, \hat y)$, it follows that
\begin{equation}\label{estimate}
\begin{array}{ll}
g(\hat y)- f(\hat x) & \leq     | \nabla_x \varphi |^\alpha F(\hat x, \nabla_x \varphi,  X_\epsilon)                               
   -| \nabla_y \varphi |^\alpha F(\hat y, -\nabla_y \varphi,  -Y_\epsilon)\\[2ex]
   & =    | \nabla_x \varphi |^\alpha \left[ F(\hat x, \nabla_x \varphi ,  X_\epsilon)  - F(\hat x, \nabla_x \varphi,  -Y_\epsilon)\right] \\[2ex]
  & \quad  +    | \nabla_x \varphi |^\alpha   \left[  F(\hat x, \nabla_x \varphi,  -Y_\epsilon)- F(\hat x, -\nabla_y \varphi,  -Y_\epsilon)
  \right]\\[2ex]
  & \quad + | \nabla_x \varphi |^\alpha   \left[  F(\hat x, -\nabla_y \varphi,  -Y_\epsilon)-  F(\hat y, -\nabla_y \varphi,  -Y_\epsilon)
  \right]\\[2ex]
 & \quad  +  \left(  | \nabla_x \varphi |^\alpha - | \nabla_y \varphi |^\alpha\right) F(\hat y, -\nabla_y \varphi,  -Y_\epsilon)
 \end{array}
 \end{equation}
 In order to estimate the four terms in the above right hand side, we need to use \eqref{matrix} and \eqref{Dx}.   It is convenient to write the function $\varphi$ in the form
 $$
 \varphi (x,y) = \Phi (x,y)   e^{ -L(d(x)+ d(y))}\, ,
 $$       
where
$$
\Phi (x,y)= M | x-y|^\gamma + L_1 |x-x_0|^2+ L_1 |y-x_0|^2\, .
$$
Hence, arguing as in   \cite{Pa}, one has
$$
\begin{array}{rl}
D^2\varphi = &   D^2\Phi \,  e^{ -L((d(x)+ d(y))} + \nabla \Phi \otimes \nabla \left( e^{ -L((d(x)+ d(y))}\right) \\[1ex]
& +    \nabla \left( e^{ -L((d(x)+ d(y))}\right) \otimes \nabla \Phi + \Phi \, D^2\left( e^{ -L((d(x)+ d(y))}\right)\, .
\end{array}
$$
Moreover,    we have
$$
\nabla \Phi (\hat x,\hat y)=M\gamma  | \hat  x-\hat  y |^{ \gamma-2} \left(\hat  x-\hat  y,-(\hat  x-\hat  y)\right) +2L_1(\hat  x-x_0,\hat  y-x_0)
$$
and
$$
D^2\Phi (\hat x,\hat y)  = \left( \begin{array}{cc}
\mathcal{B}      & - \mathcal{B}  \\
-\mathcal{B}  & \mathcal{B}  \end{array}\right) +2 L_1\,  I_{2N}\, ,
$$
with
$$
\mathcal{B}  = M \gamma | \hat  x-\hat  y |^{ \gamma-2} \left( I_N+ ( \gamma-2) \frac{(\hat  x-\hat y) \otimes  (\hat x-\hat y)}{| \hat x-\hat y |^2}\right)\, .
$$
Let us emphasize that the matrix $\mathcal{B}$ has the negative eigenvalue $M \gamma (\gamma-1) |\hat  x-\hat  y |^{\gamma -2}$ corresponding to the eigenvector $\hat  x-\hat y$.

Testing the right matrix inequality in \eqref{matrix} on the vector $\left(\hat x-\hat y, -(\hat x-\hat y)\right)$ yields
$$
\begin{array}{l}
\left( X_\epsilon+Y_\epsilon\right) (\hat x-\hat y)\cdot (\hat x-\hat y)\\[1ex]
\leq D^2 \varphi (\hat x,\hat y) \left( \begin{array}{c}
\hat x-\hat y\\
-(\hat x-\hat y)
\end{array}\right) \cdot \left( \begin{array}{c}
\hat x-\hat y\\ 
-(\hat x-\hat y)
\end{array}\right) +
\epsilon \left| D^2 \varphi (\hat x,\hat y) \left( \begin{array}{c}
\hat x-\hat y\\
-(\hat x-\hat y)
\end{array}\right) \right|^2\\[2ex]
\leq -C_1 M |\hat x-\hat y|^\gamma + \epsilon \, C_2 M^2  |\hat x-\hat y|^{2(\gamma-1)} 
\end{array}
$$
for $\delta$ sufficiently small, and with $C_1, C_2>0$ independent of $\delta$.

\noindent By choosing $\epsilon$ of the form
$$
\epsilon= \frac{C_1 |\hat x-\hat y|^{2-\gamma}}{2C_2 M}
$$
and dropping the subscript for notational convenience, we deduce that the minimal eigenvalue of the matrix $X+Y$ satisfies
$$
\lambda_{\rm min}(X+Y)\leq -C M |\hat x-\hat y|^{\gamma-2}\, .
$$
Furthermore, again the  matrix inequality \eqref{matrix}, having $\epsilon$ fixed as above and  tested on $2N$-vectors of the form $(\xi, \xi)$, with $\xi\in \R^N$, gives
$$
X+Y\leq C I_N\, 
$$
and, moreover,
$$
|X|+|Y|\leq C M |\hat x-\hat y|^{\gamma -2}\, .
$$
By assumption (H2) on $F$, we then obtain
$$
 F(\hat x, \nabla_x \varphi ,  X)  - F(\hat x, \nabla_x \varphi,  -Y)\leq \mathcal{M}^+_{a,A} (X+Y)\leq -C M |\hat x-\hat y|^{\gamma-2}\, ,
$$
so that, by \eqref{Dx},
$$
 | \nabla_x \varphi |^\alpha \left[ F(\hat x, \nabla_x \varphi ,  X)  - F(\hat x, \nabla_x \varphi,  -Y)\right] \leq
 - C M^{1+ \alpha} | \hat x-\hat y |^{ \gamma-2+ \alpha ( \gamma-1)}\, .
 $$
 Next, we observe that by \eqref{Dx} and \eqref{Dx+Dy}, and by assumption (H3), it follows that
 $$
 | \nabla_x \varphi |^\alpha   \left[  F(\hat x, \nabla_x \varphi, Y)- F(\hat x, -\nabla_y \varphi,  -Y)
  \right]\leq CM^\alpha |\hat x-\hat y|^{\alpha (\gamma-1)}\, ,
$$ 
as well as assumption (H4) implies that
$$
 | \nabla_x \varphi |^\alpha   \left[  F(\hat x, -\nabla_y \varphi,  -Y)-  F(\hat y, -\nabla_y \varphi,  -Y)
  \right]\leq  C M^{\alpha+1} |\hat x-\hat y|^{\gamma -2+\alpha (\gamma-1)+\theta_F}\, .
  $$
 Finally, we observe  that if $0< \alpha \leq 1$, then, again by \eqref{Dx+Dy}, 
$$
| |\nabla _x \varphi |^\alpha-|\nabla_y \varphi |^\alpha| \leq |\nabla_x \varphi +\nabla_y \varphi|^\alpha\leq C
$$   
and, therefore,
                  $$  \left(|\nabla_x \varphi |^\alpha-|\nabla_y \varphi|^\alpha\right) F(\hat y, -\nabla_y \varphi,  -Y)  \leq C M |\hat x-\hat y |^{  \gamma-2}\, .$$
On the other hand, if either $\alpha \leq 0$ or $\alpha> 1$,  by  the mean value's theorem we obtain that,   for some $\theta \in (0,1)$,
 $$
| |\nabla _x \varphi |^\alpha-|\nabla_y \varphi |^\alpha|  \leq | \alpha | |\nabla_x\varphi + \theta ( \nabla_x\varphi+\nabla_y\varphi)|^{ \alpha-1}|\nabla_x \varphi +\nabla_y \varphi|    \leq  C M ^{\alpha-1}  | \hat x-\hat y |^{ (\gamma-1)(\alpha-1)},
$$  
so that               
  $$  \left(|\nabla_x \varphi |^\alpha-|\nabla_y \varphi|^\alpha\right) F(\hat y, -\nabla_y \varphi,  -Y)  \leq   C  M^\alpha | \hat x-\hat y |^{\gamma-2+ ( \gamma-1) ( \alpha-1)}\, .$$
Combining all the obtained estimates with \eqref{estimate}, and observing that $\gamma-2+\alpha (\gamma-1)=(\alpha +1) \gamma -(\alpha +2)< -1$, we  conclude that
$$
g(\hat y)- f(\hat x)\leq  -CM^{1+ \alpha}  | \hat x-\hat y |^{\gamma-2+ ( \gamma-1) \alpha}\leq -C M^{1+\alpha}\delta^{\gamma-2+ ( \gamma-1) \alpha}=-\frac{C}{\delta^{\alpha+2}}
$$
and, by the boundedness of $f$ and $g$, we reach a contradiction for $\delta$ small enough.  
 \end{proof}
 We are now in a position to prove Theorem \ref{lipN}. 
 
\begin{proof}    We argue as in the previous proof, with a modification in the penalization terms of the function $\psi$. Thus, we consider,       for $x_0 \in B(\bar x, r)$ fixed,  the function 
         $$\psi(x,y) = u(x)-v(y)-\sup_{\overline{ \Omega\cap B(\bar x,\rho)}}  ( u-v) - (M \omega ( | x-y|) + L_1 |x-x_0|^2+ L_1 |y-x_0|^2) e^{ -L(d(x)+ d(y))},$$     
           where  
$$\omega ( s) =\left\{ \begin{array}{ll}
 s-\frac{s^{1+ \tau}}{ 2(1+\tau)} & \hbox{ if } 0\leq s\leq 1\, ,\\[1ex]
 1-\frac{1}{2(1+\tau)} &  \hbox{ if } s>1\, ,
 \end{array}\right.
 $$
 and 
$$0< \tau  <\left\{ \begin{array}{l}
 \inf \left\{ 
           \frac{\alpha }{2} ,\theta_F\right\} \ {\rm if}  \ \alpha \in (0,1], \\[1ex]
 \inf \left\{  \frac{1}{2},\theta_F\right\} \ \ {\rm if}  \ \alpha \notin (0,1].
         \end{array}\right.
           $$
We recall that $ \theta_F$ is  the Holder's exponent for $F$ with respect to  $x$ appearing in assumption (H4). We observe 
 that, for $0<s<1$,   $\frac{s}{2} \leq \omega (s) \leq s$, $1\geq \omega^\prime ( s) > \frac{1}{2} $ and 
$\omega^{ \prime \prime} ( s)= - \frac{\tau}{2} s^{ \tau-1}$,
hence it is negative large for $s$ small. 
               
It is clear that if we prove that $\psi \leq 0$ in $\Omega \cap B(\bar x, \rho)$, then we get the desired Lipschitz  estimate.  So we assume, by contradiction, that $\psi$ is positive somewhere, and we denote by  $(\hat x, \hat y)\in \overline{ \Omega \cap B(\bar x,\rho)}^2$  a maximum point of $\psi$. As in the proof of Proposition \ref{holderN}, we set 
$$\begin{array}{c}
\displaystyle S= \sup_{\overline{ \Omega \cap B(\bar x,\rho)}}u - \inf_{\overline{ \Omega\cap B(\bar x,\rho)}} v- \sup_{\overline{ \Omega\cap B(\bar x,\rho)}} (u-v)\, ,\\[3ex]
\displaystyle L_1 = \frac{4 S e^{2 L {\rm diam}(\Omega)}}{(\rho-r)^2}\, ,\quad L =  \frac{2}{R}.
 \end{array} $$
 Moreover, we choose
 $$
M =  \frac{2 S e^{2 L {\rm diam}(\Omega)}}{\delta}\, ,
 $$
 for $0<\delta<1$ to be fixed later small enough. 
 
 These choices for the constants $L, L_1$ and $M$ 
  force $( \hat x, \hat y)$ to satisfy 
                 $ 0<|\hat x-\hat y| < \delta,$ and $|\hat x-x_0| , | \hat y-x_0| <\frac {\rho-r}{2}$.
Furthermore, the definition  of $L$ and an analogous computation as  done for the H\"older's case, imply that  $\hat x,\ \hat y \in \Omega \cap B\left(\bar x, \frac{\rho+r}{2}\right)$.

We denote
$$\varphi(x,y) =  (M \omega ( | x-y|) + L_1 |x-x_0|^2+ L_1 |y-x_0|^2) e^{ -L(d(x)+ d(y))}, $$
and, arguing as in the proof of Proposition \ref{holderN}, we find, for every $\epsilon >0$,  $X_\epsilon$ and $Y_\epsilon$ in $\mathcal{S}$ such that 
              $$( \nabla _x \varphi( \hat x, \hat  y), X_\epsilon) \in \overline{J^{2,+}} u( \hat  x), \ ( -\nabla_y \varphi(\hat  x, \hat  y), -Y_\epsilon) \in \overline{J^{2,-}} v( \hat  y)\, ,$$
with $X_\epsilon$ and $Y_\epsilon$ satisfying \eqref{matrix}.

Next, we  observe that, for $\delta$ small enough, one has
\begin{equation}\label{Dxlip}
 \frac{M}{2} e^{-L(d( \hat x)+ d( \hat y))} \leq |\nabla_x \varphi( \hat x,  \hat y)| , |\nabla_y \varphi( \hat x,  \hat y)| \leq  2M   e^{-L(d(\hat x)+ d(\hat y))}.
\end{equation} 
Thus, we have in particular that $\nabla_x \varphi   ( \hat x,  \hat y), \ \nabla_y \varphi( \hat x,  \hat y)\neq 0$ and inequality \eqref{estimate} still holds true.                 

We now observe that, by Proposition \ref{holderN},   for any $\gamma <1$ there exists $c_{\frac{\rho+r}{2} , \gamma}>0$ such that  
$$ (M \omega ( |\hat x-\hat y|) + L_1 |\hat x-x_0|^2+ L_1 |\hat y-x_0|^2) e^{ -L(d(\hat x)+ d(\hat y))} \leq  c_{\frac{\rho+r}{2} , \gamma}|\hat x-\hat y|^\gamma\, .$$
Therefore, we obtain that 
\begin{equation}\label{Dx+Dylip}
\begin{array}{l}
|\nabla_x \varphi (\hat x,\hat y)+ \nabla_y \varphi (\hat x,\hat y)| \\[2ex]
 = \left\vert e^{ -L(d(\hat x)+ d(\hat  y))}
\left[ 2L_1 (\hat x-x_0) +2L_1(\hat y-x_0)\right. \right.\\[2ex]
\left. \left. \quad - L(\nabla d(\hat x)+ \nabla d(\hat y))  (M \omega ( |\hat x-\hat y|) + L_1 |\hat x-x_0|^2+ L_1 |\hat y-x_0|^2)\right]\right\vert \\[2ex]
\leq C | \hat x- \hat  y |^{\frac{ \gamma}{ 2}} \, ,
\end{array}
\end{equation}
for some $C>0$ depending on $r, \rho, R, S$ and $\gamma$, but independent of $\delta$. 

\noindent 
This implies,  for  $\alpha \in (0,1]$,
\begin{equation}\label{diffgrad1}
||\nabla_x \varphi |^\alpha - |\nabla_y \varphi |^\alpha| \leq |\nabla_x \varphi + \nabla_y \varphi |^\alpha \leq C  | \hat  x- \hat  y|^{ \frac{\gamma\alpha}{ 2}} ,
\end{equation}
whereas for either $\alpha \leq 0$ or $\alpha >1$,  by \eqref{Dxlip} and for some $\theta \in ]0,1[$, we get
\begin{equation}\label{diffgrad2}
\begin{array}{ll}
||\nabla_x \varphi |^\alpha - |\nabla_y \varphi |^\alpha| &  \leq 
|\alpha|     |\nabla_x \varphi  + \theta  (\nabla_x \varphi+\nabla_y \varphi )|^{ \alpha-1}|\nabla_x \varphi + \nabla_y \varphi |\\[2ex]
 & \leq C |\alpha| M^{ \alpha -1} | \hat x-\hat y|^{\frac{ \gamma}{ 2}}
 \end{array}
\end{equation}
Next, we argue as in the proof of Proposition \ref{holderN}, by using inequality \eqref{matrix} in order to estimate the matrices $X_\epsilon$ and $Y_\epsilon$. We observe that, in the present case, the hessian matrix of $\varphi$ evaluated at $(\hat x, \hat y)$ includes the matrix
$$
\left(\begin{array}{cc}
\mathcal{B} & -\mathcal{B}\\
-\mathcal{B} & \mathcal{B}
\end{array}\right)
$$
with
$$
\mathcal{B}= M\frac{\omega'(|\hat x-\hat y|)}{|\hat x-\hat y|}\left( I_N+ \left( \frac{\omega''(|\hat x-\hat y|)|\hat x-\hat y|}{\omega'(|\hat x-\hat y|)}-1\right) \frac{(\hat x-\hat y)\otimes (\hat x-\hat y)}{|\hat x-\hat y|^2}\right)\, ,
$$
which has the negative eigenvalue $M \omega''(|\hat x-\hat y|)=-\frac{M \tau}{2}|\hat x-\hat y|^{\tau-1}$ corresponding to the eigenvector $\hat x-\hat y$.

 By inequality \eqref{matrix}, fixing $\epsilon$ of the form $\frac{C}{M} |\hat  x-\hat  y|^{1-\tau}$ and dropping the subscript $\epsilon$, we deduce that
 $$
 \begin{array}{c}
 \lambda_{\rm min}(X+Y)\leq -C M |\hat  x-\hat  y|^{\tau -1}\, ,\\[2ex]
X+Y\leq C\, I_N\, ,\\[2ex]
|X|+|Y|\leq C \frac{M}{|\hat  x-\hat  y|}\, .
\end{array}
$$
Thus, by using also \eqref{Dxlip}, \eqref{Dx+Dylip}, \eqref{diffgrad1} and \eqref{diffgrad2}, estimate \eqref{estimate} in the present  case yields
$$
\begin{array}{l}
g(\hat y)-f(\hat x) \\[2ex]
\  \leq -C_1 M^{\alpha+1} |\hat x-\hat y|^{\tau -1} +C_2\left( M^\alpha |\hat x-\hat y|^{\frac{\gamma}{2}} +M^{\alpha+1} |\hat x-\hat y|^{\theta_F-1}+M |\hat x-\hat y|^{\frac{\gamma \alpha}{2}-1}\right)
\end{array}
$$
in the case $\alpha\in (0,1]$, and
$$
\begin{array}{l}
g(\hat y)-f(\hat x) \\[2ex]
\  \leq -C_1 M^{\alpha+1} |\hat x-\hat y|^{\tau -1} +C_2\left( M^\alpha |\hat x-\hat y|^{\frac{\gamma}{2}} +M^{\alpha+1} |\hat x-\hat y|^{\theta_F-1}+|\alpha | M^\alpha |\hat x-\hat y|^{\frac{\gamma}{2}-1}\right)
\end{array}
$$
in the case $\alpha \notin (0,1]$.

By the definition of $\tau$ we can select $\gamma \in (0,1)$ satisfying
$$
\gamma > \left\{
\begin{array}{l}
2  \frac{\tau}{\alpha} \ \hbox{ if } \alpha \in (0,1]\\[1ex]
2 \tau \ \hbox{ if } \alpha \notin (0,1]
\end{array}\right.
$$
so that, in any case, we obtain
$$
g(\hat y)-f(\hat x) \leq -C M^{\alpha+1} |\hat x-\hat y|^{\tau -1}\leq - \frac{C}{\delta^{2+\alpha-\tau}}\, .
$$
Since $2+\alpha>1>\tau$, the above inequality yields a contradiction to the boundedness of $g$ and $f$ for $\delta$ small enough. \end{proof}

                    \end{document}